\documentclass[12pt]{article}
\usepackage{amssymb,amsmath,amscd}
\usepackage[frenchb]{babel} 
\begin{document}

\begin{center}
{\Large \bf Foncteurs de Mackey \`a r\'eciprocit\'e}\\ 
version pr\'eparatoire (1991?)
\end{center}

\begin{center}
Bruno Kahn 
\end{center}

\newenvironment{dem}{{\bf D\'emonstration}.}{\\[5mm]}
\newenvironment{Th}{\vskip1em\noindent{\bf Th\'eor\`eme.}
\ignorespaces\it}{\par\vskip1em} 
\newenvironment{rque}{{\bf Remarque}.}{\\[5mm]} 
\newenvironment{rques}{{\bf Remarques}.\begin{enumerate}}{\end{enumerate}}
\newenvironment{ex}{{\bf Exemple}.}{\\[5mm]}
\newenvironment{exs}{{\bf Exemples}.\begin{enumerate}}{\end{enumerate}}
\newenvironment{ack}{{\bf Acknowledgements}.}{\\[5mm]}
\newenvironment{nota}{{\bf Notations}.}{\\[5mm]}
\newenvironment{qn}{{\bf Question}.}{\\[5mm]} 

\newtheorem{prop}{Proposition}[section]
\newtheorem{lemme}{Lemme}[section]
\newtheorem{defn}{D\'efinition}[section]
\newtheorem{thm}{Th\'eor\`eme}[section]
\newtheorem{cor}{Corollaire}[section]
\newtheorem{conj}{Conjecture}
\newcommand{\qq}{{\bf Q}}
\newcommand{\zz}{\bf Z\rm}
\newcommand{\zq}{{\bf Z}_{qfh}}
\newcommand{\nn}{\bf N}
\newcommand{\nq}{{\bf N}_{qfh}}
\newcommand{\oo}{\mathop\otimes}
\newcommand{\uu}{\underline}
\newcommand{\ih}{\uu{Hom}}
\newcommand{\af}{{\bf A}^1}
\newcommand{\Br}{\mbox{\rm Br }}
\newcommand{\Pic}{\mathop{\bf Pic }\nolimits}
\newcommand{\Div}{\mathop{\bf Div }\nolimits}
\newcommand{\Cor}{\mathop{\rm Cor }\nolimits}
\newcommand{\Res}{\mathop{\rm Res }\nolimits}
\newcommand{\Spec}{\mbox{\rm Spec }}
\newcommand{\et}{\mathop{\rm \acute{e}t}}
\renewcommand{\refname}{R\'ef\'erences}
\newcommand{\bH}{\mathbb{H}}
\newcommand{\fm}{\frak{m}}
\def\lhar#1#2{\smash{\mathop{\hbox to 12mm{\leftarrowfill}}
\limits^{\scriptstyle#1}_{\scriptstyle#2}}}
\def\rhar#1#2{\smash{\mathop{\hbox to 12mm{\rightarrowfill}}
\limits^{\scriptstyle#1}_{\scriptstyle#2}}}
\def\dvar#1#2{\llap{$\scriptstyle #1$}\left\downarrow
\vbox to 6mm{}\right.\rlap{$\scriptstyle #2$}}
\def\uvar#1#2{\llap{$\scriptstyle #1$}\left\uparrow
\vbox to 6mm{}\right.\rlap{$\scriptstyle #2$}}
%
%

On fixe un corps de base $k$. On consid\`ere des foncteurs de Mackey
d\'efinis sur la cat\'egorie des $k$-sch\'emas affines (donc des
$k$-alg\`ebres commutatives) par rapport aux morphismes finis et plats,
et commutant aux   limites inductives filtrantes ({\it i.e.}
``continus").

	On dit qu'un foncteur de Mackey $A$ est:

\begin{itemize}
\item {\sl cohomologique} si, pour tout morphisme fini $f$ de
degr\'e constant $d$, $A_*(f)\circ A^*(f)$ est la multiplication par
$d$;
\item {\sl additif} si, pour tous $k$-sch\'emas affines $X, X'$,
$(\iota^*,{\iota'}^*):A(X)\oplus
A(X')\stackrel{\sim}{\rightarrow}A(X\coprod X')$, o\`u $\iota$ et
$\iota'$ sont les inclusions $X\hookrightarrow X\coprod X'$ et
$X'\hookrightarrow X\coprod X'$;
\item {\sl faiblement additif} si, pour tout $k$-sch\'ema affine $S$ et
tous morphismes finis et plats $f:X\to S$ et $f':X'\to S$, on a (avec les
notations ci-dessus) $A_*(f\coprod f') = A_*(f)\circ
A^*(\iota)+A_*(f')\circ A^*(\iota')$;\\ - {\sl topologiquement invariant}
si, pour tout morphisme radiciel $f$, $A^*(f)$ est un isomorphisme;
\item (si $A$ est cohomologique:) {\sl faiblement topologiquement
invariant} si, pour toute extension $F$ de $k$ et toute $F$-alg\`ebre
locale artinienne $R$ de longueur $d$, de corps r\'esiduel $l$ de degr\'e
fini sur $F$, on a $A_*(p_R) = dA_*(p_F)\circ A^*(\iota)$, o\`u $p_R$
({\it resp.} $p_l$) est la projection
$\Spec R\to\Spec F$ ({\it resp.} $\Spec l\to\Spec F$) et $\iota$ est
l'immersion ferm\'ee $\Spec l\to\Spec R$.
\end{itemize}

Il est clair qu'un foncteur de Mackey additif est faiblement additif; de
m\^eme:

\begin{lemme}
\label{l0.1}
Un foncteur de Mackey cohomologique topologiquement invariant est
faiblement topologiquement invariant.
\end{lemme}

\begin{dem}
Avec les notations ci-dessus, supposons d'abord $l/F$ radiciel. Par
hypoth\`ese, $A^*(p_R)$ est un isomorphisme. La formule \`a d\'emontrer
est donc \'equivalente \`a
$$A_*(p_R)\circ A^*(p_R) = dA_*(p_l)\circ A^*(\iota)\circ A^*(p_R),$$
qui est claire.

Supposons maintenant $l/F$ quelconque. Soit $l_0$ la fermeture
s\'eparable de
$F$ dans $l$. Par le lemme de Hensel, $l_0$ se rel\`eve de mani\`ere
unique dans $R$ au-dessus de $F$; autrement dit, il existe un unique
$s:\Spec R\to\Spec l_0$ tel que le diagramme
$$\begin{array}{ccc}
\Spec l&\rhar{\iota}{}&\Spec R\\
\dvar{p'}{}&s\swarrow&\dvar{p_R}{}\\
\Spec l_0&\rhar{p_{l_0}}{}&\Spec F\\
\end{array}$$
soit commutatif. On a alors:
$$A_*(p_R)=A_*(p_{l_0})\circ A_*(s)$$
et
$$A_*(s)=dA_*(p')\circ A^*(\iota)$$
d'apr\`es le cas radiciel.\hfill $\Box$
\end{dem}

Toutes les courbes sur $k$ sont lisses, mais pas n\'ecessairement
compl\`etes ou irr\'eductibles. Si $U$ est un ouvert dense d'une courbe
compl\`ete $X$, on consid\'erera toujours $\Div(U)$ comme le groupe des
diviseurs de
$X$ \`a support dans $U$: cela munit $\Div$ d'une structure de
foncteur covariant pour les immersions ouvertes.

\section{Un accouplement ``divisoriel"}

Soient $A$ un foncteur de Mackey et $U$ une courbe affine sur $k$. Pour
tout point ferm\'e $x$ de $U$, on d\'efinit l'{\sl \'evaluation en $x$},
$a\mapsto a(x)$, comme le compos\'e
$$A(U)\rhar{\iota_x^*}{}A(k(x))\rhar{\Cor_{k(x)/k}}{}A(k)$$
o\`u $\iota_x$ est l'inclusion $x\mapsto U$. On l'\'etend par
lin\'earit\'e en un accouplement:
\begin{equation}
\label{acc}
\begin{array}{ccc}
A(U)\times \Div(U)&\to& A(k)\\
(a,D)&\mapsto&a(D).\\
\end{array}
\end{equation}

On a le lemme trivial suivant:

\begin{lemme}
\label{l1.0}
Soient $U'$ un ouvert dense de $U$, $a \in A(U)$, $a'$ son image
dans $A(U')$ et $D$ un diviseur de $X$ \`a support dans $U'$. Alors on a
$a(D) = a'(D)$.\hfill $\Box$
\end{lemme}

Soit $D$ un diviseur effectif de $U$, vu comme sous-sch\'ema ferm\'e de
$U$. Soit
$\iota_D$ l'inclusion de $D$ dans $U$. On a un morphisme de fonctorialit\'e
$\iota_D^*:A(U)\to A(D)$ et (puisque $D$ est fini sur $\Spec k$) un
transfert
$$A(D)\to A(k)$$
not\'e abusivement $\Cor_{D/k}$.

\begin{lemme}
\label{l1.1}
Si $A$ est cohomologique, faiblement additif et faiblement
to\-po\-lo\-gi\-que\-ment
invariant, on a pour tout diviseur effectif $D$ et tout $a \in A(U)$:
$$a(D) =\Cor_{D/k}\iota_D^*(a).$$
\end{lemme}

Cela r\'esulte du lemme \ref{l0.1}.

\begin{lemme}
\label{l1.2}
Soient $U$ et $V$ deux courbes affines sur $k$ et $f:V\to U$ un
morphisme fini de degr\'e $n$. Supposons $A$ cohomologique. Alors:\\
a) Pour tout $a \in A(U)$ et tout $D \in \Div(U)$, $(f^*a)(f^*D) =
na(D)$.\\ 
b) Pour tout $a \in A(U)$ et tout $D \in \Div(V)$, $a(f_*D) =
(f^*a)(D)$.\\ 
c) Si de plus $A$ est faiblement additif et faiblement
topologiquement invariant, on a pour tout $a \in A(V)$ et tout $D \in
\Div(U)$, $a(f^*D) = (f_*a)(D)$.
\end{lemme}

\begin{dem} Il suffit de d\'emontrer a), b) et c) lorsque $D$ est un
point ferm\'e $x$. Dans le cas a), on \'ecrit $f^*x = \sum e_iy_i$, d'o\`u
$$\begin{array}{ccl}
(f^*a)(f^*x)   &=& \sum e_i(f^*a)(y_i)\\
               &=& \sum e_i\Cor_{k(y_i)/k}\iota_{y_i}^*(f^*a)\\
               &=& \sum e_i\Cor_{k(y_i)/k}(f\circ \iota_{y_i})^*a\\
               &=& \sum e_i\Cor_{k(y_i)/k}\Res_{k(y_i)/k(x)}\iota_x^*a\\
               &=& \sum e_i[k(y_i):k(x)]\Cor_{k(x)/k}\iota_x^*a\\
               &=& n\Cor_{k(x)/k}\iota_x^*a\\
               &=& na(x).\\
\end{array}$$

(On a utilis\'e la formule $\sum e_i[k(y_i):k(x)] = n$.)\\

Dans le cas b), on a $f_*(x) = [k(x):k(y)]y$, o\`u $y = f(x)$, et
$$\begin{array}{ccl}
(f^*a)(x) &=& \Cor_{k(x)/k}\iota_x^*f^*a\\
          &=& \Cor_{k(x)/k}(f\circ \iota_x)^*a\\ 
          &=& \Cor_{k(x)/k}\iota_y^*a\\
          &=& [k(x):k(y)]\Cor_{k(y)/k}\iota_y^*a\\
          &=& a(f_*(x)),\\
\end{array}$$
puisque $A$ est cohomologique.

Pour c), notons $\Delta$ le diviseur effectif $f^*x$, $\iota_\Delta$
l'inclusion de $\Delta$ dans $U$ et $f'$ la restriction de $f$ \`a
$\Delta$. On a alors:\\
$$\iota_x^*f_*a = f'_*\iota_\Delta^*a$$
donc
$$\begin{array}{ccl}
(f_*a)(x) &=& \Cor_{k(x)/k}\iota_x^*f_*a\\
          &=&\Cor_{k(x)/k}f'_*\iota_\Delta^*(a)\\
          &=&\Cor_{\Delta/k}\iota_\Delta^*(a)\\
          &=& a(\Delta)\\
\end{array}$$
d'apr\`es le lemme \ref{l1.1}.
\end{dem}

\begin{rque} Le lemme \ref{l1.2} peut s'interpr\'eter de la mani\`ere
suivante: sous les hypoth\`eses de b), l'accouplement (i) se prolonge en
un homomorphisme
$$ev:A\oo\limits^M\Div\to A(k)$$
du foncteur de Mackey $A\oo\limits^M\Div$ vers le foncteur de Mackey
constant de valeur $A(k)$.
\end{rque}

\section{Topologie modulaire}

Soient $U$ une courbe affine (lisse) sur $k$, $X$ sa compl\'et\'ee et
$Z$ le ferm\'e compl\'ementaire. Un {\sl module} sur $U$ est un diviseur
effectif de
$X$, de support $Z$. Si $\fm$ est un module sur $U$, on lui associe un
sous-groupe $\Div^{\fm}(U)$ de $\Div(U)$:
$$\Div^{\fm}(U) = \{(f) | f \in k(X)^*; \mbox{ pour tout } x \in Z,
v_x(f-1)
\ge v_x(D)\}.$$

Les $\Div^{\fm}(U)$ forment une base de voisinages de $0$ pour une
topologie
sur $\Div(U)$: la {\sl topologie modulaire}. Si ${\fm} =
\sum\limits_{x\in Z}x$, le quotient $\Div(U)/\Div^m(U)$ s'identifie au
{\sl groupe de Picard relatif} $\Pic(X,Z)$.

Si $x$ est un point ferm\'e de $X$, on note $\hat {\cal O}_x$ le
compl\'et\'e de
${\cal O}_{X,x}$ et $\hat K_x$ le corps des fractions de $\hat {\cal
O}_x$.

\begin{defn}
\label{d2.1}
Le {\sl groupe des id\`eles} ${\cal I}(X)$ de $X$ est le produit
direct restreint des $\hat K_x^*$ relativement aux $\hat O_x^*$. Le {\sl
groupe des classes d'id\`eles} ${\cal C}(X)$ de $X$ est le quotient de
${\cal I}(X)$ par le sous-groupe image de $k(X)^*$ par le plongement
diagonal.
\end{defn}

Le lemme suivant est bien connu:

\begin{lemme}
\label{l2.1}
On a un isomorphisme canonique
$${\cal C}(X) \simeq \lim_{\longleftarrow}\Div(U)/\Div^m(U),$$
o\`u $U$ parcourt les ouverts affines de $X$ et $\fm$ parcourt les modules
sur
$U$.
\end{lemme}

\section{Foncteurs de Mackey \`a r\'eciprocit\'e}

\begin{defn}
Soient $A$ un foncteur de Mackey, $u$ une courbe affine (lisse) sur $k$,
$a\in A(U)$ et $\fm$ un module sur $U$. On dit que $a$ admet le module
$\fm$ si la restriction de l'application $\Div(U)\to A(F)$ induite par $a$
\`a
$\Div^{\fm}(U)$ est identiquement nulle.
\end{defn}

\begin{defn}
\label{d3.1}
Soient $A$ un foncteur de Mackey et $U$ une courbe affine (lisse)
sur $k$. On dit que $A$ {\sl v\'erifie la $k$-r\'eciprocit\'e sur $U$} si
l'accouplement {\rm (\ref{acc})} est continu pour la topologie modulaire
sur
$\Div(U)$ et la topologie discr\`ete sur $A(U)$ et $A(k)$. On dit que
$A$ {\sl v\'erifie la $k$-r\'eciprocit\'e forte sur $U$} si, de plus,
l'accouplement
$A(U)\times\Div(U)\to A(k)$ restreint \`a $A(U)\times\Div^{\fm}(U)$
est identiquement nul pour ${\fm} = \sum\limits_{x\in Z}x$.\\
On dit que $A$ {\sl v\'erifie la $k$-r\'eciprocit\'e} (resp. {\sl
v\'erifie la $k$-r\'eciprocit\'e forte}) s'il la v\'erifie sur toute
courbe
$U$ sur $k$.\\
On dit que $A$ {\sl v\'erifie la r\'eciprocit\'e} (resp. {\sl
v\'erifie la r\'eciprocit\'e forte}) s'il v\'erifie la
$l$-r\'eciprocit\'e (resp. la $l$-r\'eciprocit\'e forte) pour toute
extension finie $l$ de $k$.
\end{defn}

\begin{rque}
Supposons $k$ alg\'ebriquement clos. Dans la terminologie de
\cite[ch. III]{Se}, la d\'efinition de la $k$-r\'eciprocit\'e signifie que, pour
tout $a \in A(U)$, l'application $x\mapsto a(x)$ de $U$ dans $A(k)$
poss\`ede un module. \end{rque}

\begin{lemme}
\label{l3.1}
Supposons $A$ cohomologique, faiblement additif et faiblement
topologiquement invariant. Pour que $A$ v\'erifie la $k$-r\'eciprocit\'e
forte, il faut et il suffit qu'il soit $k$-invariant par homotopie, i.e.
que
$A(k)\stackrel{\sim}{\to}A({\bf A}^1_k)$.
\end{lemme}

En effet, supposons que $A$ v\'erifie la $k$-r\'eciprocit\'e forte.
Prenons $U = {\bf A}^1_k$. On a alors $X = {\bf P}^1_k$ et $\Pic({\bf
P}^1_k,\{\infty\}) = {\bf Z}$. Soient $a \in A(U)$ et $a_0 = a(0)$: on
peut voir $a_0$ comme un
\'el\'ement de $A(U)$ via l'homomorphisme $A(k)\to A(U)$ d\'eduit du
morphisme structural. Montrons que $a = a_0$: cela montrera que $A(k)\to
A(U)$ est surjectif, donc bijectif puisque c'est {\it a priori} une
injection scind\'ee. Quitte \`a remplacer $a$ par $a-a_0$ on peut
supposer $a_0 = 0$. Soit $x$ un point ferm\'e de $U$, de degr\'e $d$: il
existe une (unique) fonction rationnelle $f \in k(U)^*$ telle que $(f) =
x-d0$ et que $f(\infty) = 1$ (c'est
$f(t) = P(t)/t^d$, o\`u $P$ est le polyn\^ome minimal de $x$). Par
hypoth\`ese, on a $0 = a((f)) = a(x)-da(0) = a(x)$.

Supposons maintenant $A$ invariant par homotopie. Soient $U$ une courbe
affine et irr\'eductible sur $k$, $X$ sa compl\'et\'ee, $a \in A(U)$ et
$f \in k(X)^*$ une fonction rationnelle (suppos\'ee non constante), telle
que $f(x) = 1$ pour tout
$x \notin U$. Consid\'erons $f$ comme un morphisme (fini et plat) de $X$
dans
${\bf P}^1_k$, donc comme un morphisme fini et plat d'un ouvert $U'$ de
$U$ dans
${\bf P}^1_k-\{1\}$. On a alors $(f) = f^*(0-\infty)$, donc d'apr\`es
les lemmes
\ref{l1.0} et \ref{l1.2}:
$$a((f)) = a'((f)) = a'(f^*(0-\infty)) = (f_*a')(0-\infty),$$
o\`u $a'$ est l'image de $a$ dans $U'$. En appliquant l'invariance par
homotopie
\`a ${\bf P}^1_k-\{1\} \simeq {\bf A}^1_k$, on trouve $(f_*a')(0-\infty)
= 0$, donc $a((f)) = 0$.\\

{\bf Question}. Est-il vrai que $A$ v\'erifie la $k$-r\'eciprocit\'e si
et seulement si il la v\'erifie sur les ouverts de ${\bf P}^1_k$? ({\it
c.f.}
\cite[ch III, prop. 9]{Se}).

\begin{defn}
On dit qu'un foncteur de Mackey $A$ {\sl v\'erifie la r\'eciprocit\'e}
(resp. {\sl la r\'eciprocit\'e forte}) s'il v\'erifie la
$l$-r\'eciprocit\'e (resp. la
$l$-r\'eciprocit\'e forte) pour toute extension finie $l$ de $k$.
\end{defn}

\section{R\'eciprocit\'e et symboles locaux}

\begin{defn}
\label{d4.1}
Soit $A$ un foncteur de Mackey. Un {\sl symbole local} associ\'e \`a $A$
est la donn\'ee $\partial$, pour toute $k$-alg\`ebre de valuation
discr\`ete d'origine g\'eom\'etrique $\cal O$, de corps r\'esiduel $l$
alg\'ebrique sur $k$ et de corps des fractions $K$, d'un accouplement
$\partial_{\cal O}:A(K)\times K^*\to A(l)$ ayant les propri\'et\'es
suivantes:\\ (i) $\partial_{\cal O}$ est continu pour la topologie
naturelle de $K^*$ et pour les topologies discr\`etes de $A(K)$ et
$A(l)$.\\ (ii) Soient $a$ un \'el\'ement de $A({\cal O})$, $a'$ son image
dans $A(K)$ et
$\overline a$ son image dans $A(l)$. Alors, pour tout $x \in K^*$, on a 
$\partial_{\cal O}(a',x) = v(x)\overline a$, o\`u $v$ est la valuation
associ\'ee
\`a $\cal O$.\\
(iii) Soit ${\cal O}'$ une extension finie, int\'egralement close de
$\cal O$: c'est un anneau principal semi-local. Soient $K'$ son corps des
fractions,
${\cal O}_i$ les localis\'es de ${\cal O}'$ en ses id\'eaux maximaux,
$l_i$ le corps r\'esiduel de ${\cal O}_i$, $e_i$ l'indice de ramification
de ${\cal O}_i/{\cal O}$. Alors:\\
(iii1) Pour $(a,x) \in A(K)\times K^*$ et pour tout $i$,
$$\partial_{{\cal O}_i}(\Res_{K'/K}a,x) = e_i\Res_{l_i/l}\partial_{\cal
O}(a,x);$$
(iii2) Pour $(a,x) \in A(K)\times {K'}^*$, 
$$\partial_{\cal O}(a,N_{K'/K}x) =\sum\Cor_{l_i/l}\partial_{{\cal
O}_i}(\Res_{K'/K}a,x);$$
(iii3) Pour $(a,x) \in A(K')\times K^*$, 
$$\partial_{\cal O}(\Cor_{K'/K}a,x)
=\sum\Cor_{l_i/l}\partial_{{\cal O}_i}(a,x).$$
On dit que $\partial$ est un {\sl symbole local fort} s'il v\'erifie la
condition suivante (qui implique (i)):\\
(i) fort: $\partial_{\cal O}$ est nul sur $A(K)\times U_1$, o\`u $U_1$
d\'esigne le groupe des unit\'es principales de $\cal O$.
\end{defn}

\begin{rques}
\item L'existence d'un symbole local impose des restrictions sur la
structure de groupe ab\'elien des valeurs de $A$. Par exemple, la
propri\'et\'e (ii) implique que l'application naturelle $A({\cal O})\to
A(l)$ peut se prolonger
\`a $A(K)$.
\item Les conditions (iii) impliquent que $\partial_{\cal O}$ se
prolonge en un homomorphisme (encore not\'e) $\partial_{\cal
O}:A\oo\limits^M {\bf G}_m(K)\to A(l)$, ayant les propri\'et\'es
suivantes (sous les hypoth\`eses de (iii)):\\
(iii1) bis: Pour tout $b \in A\oo\limits^M {\bf G}_m(K)$ et pour tout
$i$,
$\partial_{{\cal O}_i}(\Res_{K'/K}b) = e_i\Res_{l_i/l}\partial_{\cal
O}(b)$;\\ 
(iii2) bis: Pour tout $b \in A\oo\limits^M {\bf G}_m(K')$,
$\partial_{\cal O} (\Cor_{K'/K}b) = \sum\Cor_{l_i/l}\partial_{{\cal
O}_i}(b)$.
\end{rques}

\begin{lemme}
\label{l4.1}
Pour se donner un symbole local (resp. un symbole local fort), il
suffit de se donner une famille d'accouplements $\partial_{\cal O}$
ayant les propri\'et\'es (i)-(iii) (resp. (i) fort-(iii)) pour $\cal O$
parcourant les
$k$-alg\`ebres de valuation discr\`ete hens\'eliennes [d'origine
g\'eom\'etrique] dont le corps r\'esiduel est fini sur $k$.
\end{lemme}

Supposons donn\'ee une telle famille. Soit $\cal O$ une $k$-alg\`ebre de
valuation discr\`ete d'origine g\'eom\'etrique, de corps r\'esiduel $l$
et de corps de fractions $K$. Soient ${\cal O}^h$ la hens\'elis\'ee de
$\cal O$ (de corps r\'esiduel $l$) et $K^h$ le corps des fractions de
$O^h$. On d\'efinit
$\partial_{\cal O}$ comme le compos\'e:
$$A(K)\times K^*\to A(K^h)\times K^{h*}\to A(l),$$
l'application de droite \'etant $\partial_{{\cal O}^h}$. On v\'erifie sans
peine que les propri\'et\'es (i)--(iii) sont v\'erifi\'ees.

\begin{thm}
\label{t4.1}
Soit $A$ un foncteur de Mackey cohomologique, faiblement additif et
faiblement topologiquement invariant. Les conditions suivantes sont
\'equivalentes:\\
a) $A$ v\'erifie la r\'eciprocit\'e (resp. la r\'eciprocit\'e
forte).\\
b) $A$ poss\`ede un symbole local (resp. un symbole local fort)
$\partial$, v\'erifiant la condition suivante: pour toute extension finie
$l$ de
$k$, toute courbe $X$ lisse, compl\`ete, irr\'eductible sur $l$, de corps
des fonctions
$K$, et tout $(a,f) \in A(K)\times K^*$, on a 
$$\sum \Cor_{l(x)/l}\partial_x(a,f) = 0,$$
o\`u $x$ parcourt les points ferm\'es de $X$ et, pour tout $x$,
$\partial_x$ est le symbole local associ\'e \`a ${\cal O}_{X,x}$. De
plus, sous ces conditions, le symbole local $\partial$ est unique.
\end{thm}

\begin{dem} Montrons que b) $\Rightarrow$ a). Soient $U$ un ouvert
affine non vide de $X$, $a \in A(U)$ et $a'$ l'image de $a$ dans $A(K)$.
Pour tout $x \in U$ et tout $f \in K^*$, on a (propri\'et\'e (ii) de la
d\'efinition \ref{d4.1}):
$$\partial_x(a',f) = v_x(f)\iota_x^*a.$$

En particulier, supposons que $(f)$ soit \`a support dans $U$. Alors
$$\sum_{x\in U}\Cor_{k(x)/k}\partial_x(a,f) = \sum_{x\in
U}v_x(f)\Cor_{k(x)/k}\iota_x^*a = a((f)).$$

Par la propri\'et\'e (i) de la d\'efinition \ref{d4.1}, il existe un module
$\fm$ pour $U$ tel que, si $(f) \in \Div^{\fm}(U)$, on ait
$\partial_x(a,f) = 0$ pour tout $x \notin U$. On a alors:
$$a((f)) = \sum_{x\in U}\Cor_{k(x)/k}\partial_x(a,f) = \sum_{x\in
X}\Cor_{k(x)/k}\partial_x(a,f) = 0,$$
donc l'accouplement (\ref{acc}) est continu pour la topologie modulaire. Si
$\partial$ est un symbole local fort, on peut choisir ci-dessus ${\fm} = 
\sum\limits_{x\notin U}x$, donc $A$ v\'erifie la $k$-r\'eciprocit\'e forte.\\

Montrons que a) $\Rightarrow$ b). Comme $A$ est continu, on a $A(K) =
\lim\limits_{\longrightarrow}A(U)$, o\`u $U$ d\'ecrit les ouverts non
vides de $X$.
Gr\^ace au lemme \ref{l2.1}, on obtient donc un accouplement continu:
$$A(K)\times {\cal C}(X)\to A(k).$$

Soit $x$ un point ferm\'e de $X$. Par restriction \`a (l'image de) $\hat
K_x^*$, on obtient un accouplement local continu:
$$\delta_x:A(K)\times\hat K_x^*\to A(k),$$
tel que $\delta_x(a,f) = v_x(f)a(x)$ si $a$ ``provient de ${\cal
O}_{X,x}$". D'o\`u encore par restriction un autre accouplement local
continu
$$A(K)\times K_x^{h*}\to A(k),$$
o\`u $K_x^h$ est le hens\'elis\'e de $K$ en $x$.\\

Choisissons $X = {\bf P}^1_k$, $x = 0$. Alors $\hat K_x$ est le corps des
s\'eries formelles $k((t))$, et $K_x^h$ est le sous-corps $k\{\{t\}\}$
des s\'eries formelles alg\'ebriques sur $k$, corps des fractions du
hens\'elis\'e
$k\ll t\gg$ de $k[t]$ en $0$. Soit $Y\to X$ un rev\^etement non
ramifi\'e et d\'ecompos\'e en $x$, et soit $y$ un point de $Y$ au-dessus
de $X$, de corps r\'esiduel $k$. Si $L$ est le corps des fonctions de
$Y$, on a
$K_x^h\stackrel{\sim}{\to}L_y^h$, et le lemme \ref{l1.2} b) montre que le
diagramme
$$\begin{array}{ccccc}
A(L)&\times& L_y^{h*}&	\to&	A(k)\\
&\uparrow & & &		\mid\mid\\
A(K)&\times& K_x^{h*}&	\to&	A(k)\\
\end{array}$$
est commutatif. En passant \`a la limite, on obtient un accouplement
continu
$$\partial_{k\ll t\gg}:A(k\{\{t\}\})\times k\{\{t\}\}^*\to A(k),$$
ayant la propri\'et\'e (ii) de la d\'ef. \ref{d4.1} (observer que
$k\{\{t\}\} =
\lim\limits_{\longrightarrow} k(Y)$, o\`u $Y$ parcourt les rev\^etements
de ${\bf
P}^1_k$ du type ci-dessus).\\

En r\'ep\'etant cette op\'eration avec pour base une extension finie
arbitraire
$l$ de $k$ et en tenant compte du lemme \ref{l4.1}, on obtient un symbole
local: en effet, toute $k$-alg\`ebre de valuation discr\`ete hens\'elienne,
de corps r\'esiduel $l$ est de la forme $l\ll t\gg$ (on v\'erifie
facilement les propri\'et\'es (iii) de la d\'efinition \ref{d4.1} \`a
l'aide du lemme \ref{l1.2}). Reste \`a v\'erifier la formule du th.
\ref{l4.1} b). On peut supposer $l=k$. Avec les notations ci-dessus,
$(a,f) \in A(K)\times K^*$ et $x
\in X$, on a:
$$\delta_x(a,f) = \Cor_{k(x)/k}\partial_x(a,f).$$

On en d\'eduit:
$$\sum_{x\in X}\Cor_{k(x)/k}\partial_x(a,f) = \sum_{x\in X}\delta_x(a,f)
= 0,$$ 
puisque la classe de $f$ est triviale dans ${\cal C}(X)$.\\

Enfin, l'unicit\'e du symbole $\partial$ r\'esulte de la d\'emonstration
de b)
$\Rightarrow$ a).
\end{dem}

\section{Exemples de foncteurs \`a r\'eciprocit\'e}

\begin{defn}
\label{d5.1}
Un foncteur de Mackey $A$ est {\sl propre} si, pour toute $k$-alg\`ebre
de valuation discr\`ete $\cal O$, de corps des fractions $K$, $A({\cal
O})\to A(K)$ est surjective.
\end{defn}

\begin{prop}
\label{p5.1} Soit $A$ un foncteur de Mackey propre. Pour que $A$ puisse
\^etre muni d'un symbole local $\partial$, il faut et il suffit que, pour
toute
$k$-alg\`ebre de valuation discr\`ete hens\'elienne $\cal O$, de corps des
fractions $K$ et de corps r\'esiduel $l$, l'application
$\mbox{\rm Ker }(A({\cal O})\to A(K))\to A(l)$ soit identiquement nulle.
Le symbole
$\partial$ est alors unique, et donn\'e par $\partial_{\cal O}(a,x) =
v(x)\overline a$ (notations de la d\'ef. \ref{d4.1}); il est fort.
\end{prop}

Cela r\'esulte imm\'ediatement de la propri\'et\'e (ii) de la d\'ef.
\ref{d4.1}.

\begin{cor}
Soit $A$ un foncteur de Mackey propre, cohomologique, faiblement
additif et faiblement topologiquement invariant. Pour que $A$ v\'erifie la
r\'eciprocit\'e, il faut et il suffit qu'il soit invariant par homotopie.
\end{cor}

Cela r\'esulte de la prop. \ref{p5.1} et du lemme \ref{l3.1}.\\

Les deux lemmes suivants ne pr\'esentent aucune difficult\'e.

\begin{lemme}
\label{l5.1}
Soient $A$ et $B$ deux foncteurs de Mackey.\\
a) Si $A$ et $B$ v\'erifient la r\'eciprocit\'e (resp. la
r\'eciprocit\'e forte), il en est de m\^eme pour $A\oplus B$.\\
b) Soit $\varphi:A\to B$ un morphisme de foncteurs
de Mackey.\\
b1) Si $\varphi$ est surjectif et si $A$ v\'erifie la r\'eciprocit\'e
(resp\dots), il en est de m\^eme pour $B$.\\
b2) Si $A(k)\to B(k)$ est injectif et si $B$ v\'erifie
la r\'eciprocit\'e (resp\dots), il en est de m\^eme pour $A$.
\end{lemme}

\begin{lemme}
\label{l5.2} Soient $(A_i)$ un syst\`eme inductif de foncteurs de Mackey
et $A = \lim\limits_{\longrightarrow}A_i$. Si les $A_i$ v\'erifient la
r\'eciprocit\'e
(resp\dots), il en est de m\^eme pour $A$.
\end{lemme}

Soit $0\to A\to B\to C\to 0$ une suite exacte de foncteurs de Mackey. Si
$A$ et
$C$ v\'erifient la r\'eciprocit\'e, j'ignore s'il en est de m\^eme en
g\'en\'eral pour $B$. C'est cependant le cas si $A$ et $C$ v\'erifient la
r\'eciprocit\'e forte:

\begin{prop}
\label{p5.3} Soit $0\to A\to B\to C\to 0$ une suite exacte de foncteurs de
Mackey cohomologiques, faiblement additifs et faiblement topologiquement
invariants. Si $A$ et $C$ v\'erifient la $k$-r\'eciprocit\'e forte, il
en est de m\^eme pour $B$.
\end{prop}
	
En effet, d'apr\`es le lemme \ref{l3.1}, $A$ et $C$ sont invariants par
homotopie. Le diagramme commutatif aux lignes exactes
$$\begin{array}{ccccccccc}
0	&\to& 	A(k)									  &\to& 	B(k)									  &\to& 	C(k)	          &\to& 	0\\
			&&\downarrow&&		\downarrow&&		\downarrow&&\\
0	&\to& 	A({\bf A}^1_k)	&\to& 	B({\bf A}^1_k)	&\to& 	C({\bf A}^1_k)	&\to& 	0\\
\end{array}$$
montre alors que $B$ est invariant par homotopie.

\begin{thm}
\label{t5.1} {\rm (Rosenlicht)} Un foncteur de Mackey d\'efini par un
groupe alg\'ebrique commutatif (resp. par une vari\'et\'e
semi-ab\'elienne) v\'erifie la r\'e\-ci\-pro\-ci\-t\'e (resp. la
r\'eciprocit\'e forte).
\end{thm}

\begin{dem} Si $k$ est alg\'ebriquement clos, cela r\'esulte de \cite[ch.
III]{Se}. En g\'en\'eral, soient $A$ un groupe alg\'ebrique commutatif,
$U$ une courbe affine sur $k$ et $a \in A(U)$. Soient $X$ la
compl\'et\'ee de $U$ et $Z = X-U$. Soient $\overline k$ une cl\^oture
alg\'ebrique de $k$, $\overline U = U\oo_k\overline k$. et $\overline Z =
Z\oo_k\overline k$. Notons $a'$ l'image de
$a$ dans $A(\overline U)$. Il existe un module $\overline {\fm}$ de
support
$\overline Z$ tel que $a'(\Div^{\overline {\fm}}(\overline U)) = 0$. Soit 
$\overline {\fm}'$ le satur\'e de $\overline {\fm}$ pour l'action de
$\mbox{\rm Gal }(\overline k/k)$: si $p$ est l'exposant
ca\-rac\-t\'e\-ris\-ti\-que de
$k$, il existe un entier $n \ge 0$ tel que ${\fm} = p^n\overline {\fm}'$
soit rationnel sur $k$. A fortiori, on a $a'(\Div^{\fm}(\overline U)) = 0$,
d'o\`u $a(\Div^{\fm}(U)) = 0$ (lemme \ref{l1.2} b)).\\

Supposons maintenant que $A$ soit une vari\'et\'e semi-ab\'elienne,
c'est-\`a-dire une extension d'une vari\'et\'e ab\'elienne par un tore.
Alors
$A$ est invariant par homotopie, donc v\'erifie la r\'eciprocit\'e forte
d'apr\`es le lemme \ref{l3.1}.
\end{dem}

\begin{thm}
\label{t5.2}
Soit $C^\cdot$ un complexe born\'e \`a gauche de faisceaux de groupes
ab\'eliens sur le grand site \'etale de $\Spec k$. Supposons que les
faisceaux d'ho\-mo\-lo\-gie de $C^\cdot$ soient de torsion premi\`ere \`a
la caract\'eristique de $k$. Alors, pour tout $i \in {\bf Z}$, le foncteur
de Mackey
$X\mapsto {\bH}^i(X_{\et},C^\cdot)$ v\'erifie la
r\'eciprocit\'e forte.
\end{thm}

En effet, il suffit de montrer que ce foncteur de Mackey est
cohomologique, additif, topologiquement invariant et invariant par
homotopie (lemme
\ref{l3.1}). L'additivit\'e est \'evidente. Pour le reste, supposons
d'abord
$C^\cdot$ concentr\'e en degr\'e z\'ero: cela r\'esulte alors de [SGA 4
XVII (6.2.3), VIII (1.1.2), XV (2.2.2)]. Le cas g\'en\'eral r\'esulte de
celui-ci et de la suite spectrale d'hyper\-cohomologie.

\begin{thm}
\label{t5.3}
Soient $X$ une vari\'et\'e lisse sur $k$ et $i, j$ deux entiers $\ge 0$.
Notons $A$ le foncteur de Mackey d\'efini par $A(Y) = H^i((X\times_k
Y)_{Zar},{\cal K}_j)$. Alors $A$ v\'erifie la r\'eciprocit\'e forte.
\end{thm}

En effet, $A$ est cohomologique, additif, topologiquement invariant et
invariant par homotopie ([\dots]).

\begin{cor}
Avec les notations du th. \ref{t5.3}, $Y\mapsto CH^i(X\times_k Y)$
v\'erifie la r\'eciprocit\'e forte.
\end{cor}
	
C'est le cas particulier $j = i$ ([\dots]).\\

\begin{thm}
\label{t5.4}
Pour tout $i \ge 0$, le foncteur de Mackey $Y\mapsto \Omega^i_{X/{\bf Z}}$
v\'erifie la r\'eciprocit\'e (mais non la r\'eciprocit\'e forte).
\end{thm}

D\'emonstration.

\section{Produits tensoriels}

Soient $(A,\partial)$ et $(B,\partial')$ deux foncteurs de Mackey munis de
symboles locaux. On aimerait munir $A\oo\limits^M B$ d'un symbole local
$\partial'' = \partial\oo \partial'$. Si $\cal O$ est une $k$-alg\`ebre de
valuation discr\`ete, de corps des fractions $K$ et de corps r\'esiduel
$l$, et si $\pi$ est une uniformisante de $\cal O$, on interpr\`ete les
homomorphismes
$$s_\pi:A(K)\to A(l)$$
$$s'_\pi:B(K)\to B(l)$$
donn\'es par $s_\pi(a) = \partial_{\cal O}(a,\pi), s'_\pi(b) =
\partial'_{\cal O}(b,\pi)$, comme des homomorphismes de
sp\'ecialisation. On d\'esire alors d\'efinir $\partial''$ de telle sorte
que, en posant $s''_\pi(c) =
\partial''_{\cal O}(c,\pi)$, on ait identiquement (pour toute
uniformisante
$\pi$):
$$s''_\pi(a\oo b) = s_\pi(a)\oo s'_\pi(b).$$

Ceci impose des relations sur les $s_\pi(a)\oo s'_\pi(b)$, qui en
g\'en\'eral ne sont pas v\'erifi\'ees. On est donc conduit \`a quotienter
le foncteur
$A\oo\limits^M B$ par ces nouvelles relations.

\end{document}